\documentclass{article}
\usepackage{amsfonts}
\usepackage{amsmath}

\input{tcilatex}
\begin{document}

\begin{center}
\bigskip

\bigskip \textbf{THE INVERSE SURFACES OF TANGENT DEVELOPABLES WITH RESPECT
TO S}$_{{\LARGE c}}$\textbf{(r)}

\textbf{\bigskip }

\textbf{\ M. Evren AYDIN and Mahmut ERG\"{U}T}

Department of Mathematics, Firat University

Elazig, 23119, Turkey

email: aydnevren@gmail.com\\[0pt]

\bigskip
\end{center}

\textbf{Abstract. }In this paper, we define the inverse surface of a tangent
\linebreak developable surface with respect to the sphere $S_{c}\left(
r\right) $ with the center $c\in \mathbb{E}^{3}$ and the radius $r$\ in
3-dimensional Euclidean space $\mathbb{E}^{3}$. We obtain the \linebreak
curvatures, the Christoffel symbols and the shape operator of this inverse
\linebreak surface by the help of these of the tangent developable surface.
Morever, we give some necessary and sufficient conditions regarding the
inverse surface being flat and minimal.\newline

\bigskip


\textbf{Msc.} 11A25, 53A04, 53A05.\newline

\bigskip

\section{\textbf{Introduction}}

The last ten years, the developable ruled surfaces are studied by
many\linebreak\ mathematicians. Developable surfaces are a type of important
and fundamental surfaces universally used in industry design. Different
methods have been presented for the design of developable surfaces. The use
of developable surfaces in ship design is of engineering importance because
they can be easily manufactured without stretching or tearing, or without
the use of heat treatment. In some cases, a ship hull can be entirely
designed with the use of developable surfaces. See, $\cite{4,8,9,10,11}.$

On the other hand, a conformal map is a function which preserves the angles.
The conformal mapping is an important technique used in complex analysis and
has many applications in different physical situations.

An inversion with respect to the sphere $S_{c}\left( r\right) $ with the
center $c\in \mathbb{E}^{3}$ and the radius $r$ given by%
\begin{equation*}
c+\frac{r^{2}}{\left\Vert p-c\right\Vert ^{2}}\left( p-c\right) ,
\end{equation*}%
$p\in \mathbb{E}^{3},$ is a conformal mapping and also is differentiable. In 
$\mathbb{E}^{3},$ the inversion is a transformation defining between open
subsets of $\mathbb{E}^{3}.$

In this paper, we firstly tell inversions and inversions of surfaces in $%
\mathbb{E}^{3}.$ Next, we give the fundamental forms, the curvatures(Gauss
and mean), the shape operator and the Christoffel symbols of the tangent
developable. Finally, using by these properties, we obtain these of the
inverse surface of the tangent developable.

\bigskip

\section{\textbf{Basic notions of inverse surfaces}}

\bigskip

Let $c$ $\in $ $\mathbb{E}^{3}$ and $r$ $\in 
\mathbb{R}
^{+}$. We denote that $\left( \mathbb{E}^{3}\right) ^{\ast }=$ $\mathbb{E}%
^{3}-\left\{ c\right\} .$ Then, an inversion of $\mathbb{E}^{3}$ with the
center $c\in \mathbb{E}^{3}$ and the radius $r$\ is the map 
\begin{equation*}
\Phi \left[ c,r\right] :\left( \mathbb{E}^{3}\right) ^{\ast }\longrightarrow
\left( \mathbb{E}^{3}\right) ^{\ast }
\end{equation*}%
given by 
\begin{equation}
\Phi \left[ c,r\right] \left( p\right) =c+\frac{r^{2}}{\left\Vert
p-c\right\Vert ^{2}}\left( p-c\right) .  \tag{2.1}
\end{equation}

\textbf{Definition 2.1.} $\left( \cite{7}\right) $ Let $\Phi \left[ c,r%
\right] $ be an inversion with the center $c$ and the radius $r$ \ Then, the
tangent map of $\Phi $ at $p\in \left( \mathbb{E}^{3}\right) ^{\ast }$ is
the map%
\begin{equation*}
\Phi _{\ast p}:T_{p}\left( \left( \mathbb{E}^{n}\right) ^{\ast }\right)
\longrightarrow T_{\Phi \left( p\right) }\left( \left( \mathbb{E}^{n}\right)
^{\ast }\right)
\end{equation*}%
\ given by%
\begin{equation*}
\Phi _{\ast p}\left( v_{p}\right) =\frac{r^{2}v_{p}}{\left\Vert
p-c\right\Vert ^{2}}-\frac{2r^{2}\left\langle \left( p-c\right)
,v_{p}\right\rangle }{\left\Vert p-c\right\Vert ^{4}}\left( p-c\right) ,
\end{equation*}%
where $v_{p}\in T$ $_{p}\left( \left( \mathbb{E}^{3}\right) ^{\ast }\right)
. $\newline

\bigskip

Now, let us assume that $\mathfrak{X}:U\subset \mathbb{E}^{2}\longrightarrow
\left( \mathbb{E}^{3}\right) ^{\ast }$ is the parametrization of a surface.
The inverse surface of $\mathfrak{X}$ with respect to $\Phi \left[ c,r\right]
$ is the surface given by

\begin{equation}
\mathfrak{Y}=\Phi \left[ c,r\right] \circ \mathfrak{X},  \tag{2.2}
\end{equation}

Throughout this paper, we assume that $\Phi $ is an inversion of $\mathbb{E}%
^{3}$ with the center $c$ and the radius $r,$ $\mathfrak{X}$ is a patch in $%
\left( \mathbb{E}^{3}\right) ^{\ast }$and $\mathfrak{Y}$ is inverse patch of 
$\mathfrak{X}$ with respect to $\Phi $.

Let\textbf{\ }$I_{\mathfrak{X}},$ $II_{\mathfrak{X}}$ \ and $K_{\mathfrak{X}%
},$ $H_{\mathfrak{X}}$ be the first and second fundamental forms\ and the
curvatures (Gauss and mean) of $\mathfrak{X}$, and let $I_{\mathfrak{Y}}$, $%
II_{\mathfrak{Y}}$ and $K_{\mathfrak{Y}},$ $H_{\mathfrak{Y}}$ be these of $%
\mathfrak{Y},$ respectively. From $\cite{1},$ we have \textbf{\ }%
\begin{equation}
I_{\mathfrak{Y}}\circ \Phi _{\ast }=\lambda ^{2}I_{\mathfrak{X}},  \tag{2.3}
\end{equation}%
\begin{equation}
II_{\mathfrak{Y}}\circ \Phi _{\ast }=-\lambda II_{\mathfrak{X}}-2\delta I_{%
\mathfrak{X}},  \tag{2.4}
\end{equation}%
\begin{equation}
K_{\mathfrak{Y}}=\dfrac{1}{\lambda ^{2}}K_{\mathfrak{X}}+\dfrac{4}{r^{2}}%
\lambda ^{-1}\eta H_{\mathfrak{X}}+\dfrac{4}{r^{4}}\eta ^{2},  \tag{2.5}
\end{equation}%
\begin{equation}
H_{\mathfrak{Y}}=-\dfrac{1}{\lambda }H_{\mathfrak{X}}-\frac{2\eta }{r^{2}}, 
\tag{2.6}
\end{equation}%
where $\lambda =\tfrac{r^{2}}{\left\Vert \mathfrak{X}-c\right\Vert ^{2}}$, $%
\delta =\frac{2r^{2}\left\langle U_{\mathfrak{X}},\left( \mathfrak{X}%
-c\right) \right\rangle }{\left\Vert \mathfrak{X}-c\right\Vert ^{4}}$ and $%
\eta =\left\langle U_{\mathfrak{X}},\left( \mathfrak{X}-c\right)
\right\rangle $.

\bigskip

\section{\textbf{The tangent developable surface}}

\bigskip

Let $\gamma :I\subset 
\mathbb{R}
\longrightarrow \mathbb{E}^{3}$ be a curve with arc-length $s$ and $\left\{
T,N,B\right\} $ be Frenet frame along $\gamma $. Denote by $\kappa $ and $%
\tau $ the curvature and the torsion of the curve $\gamma $, respectively.
Then we have Frenet formulas%
\begin{eqnarray*}
T^{\prime }\left( s\right)  &=&\kappa \left( s\right) N\left( s\right) , \\
N^{\prime }\left( s\right)  &=&-\kappa \left( s\right) T\left( s\right)
+\tau \left( s\right) B\left( s\right) , \\
B^{\prime }\left( s\right)  &=&-\tau \left( s\right) N\left( s\right) .
\end{eqnarray*}%
The tangent developable of $\gamma $ is a ruled surface parametrized by%
\begin{equation}
\mathfrak{M}\left( s,u\right) =\gamma \left( s\right) +uT\left( s\right) , 
\tag{3.1}
\end{equation}%
where $T$ is unit tangent vector field of $\gamma .$ As it is known, the
coefficients of the first and second fundamental forms of the surface $%
\mathfrak{M}\left( s,u\right) $ have following%
\begin{equation}
E_{\mathfrak{M}}=1+\left( u\kappa \right) ^{2},\text{ }F_{\mathfrak{M}}=G_{%
\mathfrak{M}}=1,  \tag{3.2}
\end{equation}%
and%
\begin{equation}
e_{\mathfrak{M}}=-sgn\left( u\kappa \right) \left( u\kappa \tau \right) ,%
\text{\ }f_{\mathfrak{M}}=g_{\mathfrak{M}}=0.  \tag{3.3}
\end{equation}%
The normal vector field of the surface $\mathfrak{M}\left( s,u\right) $ is
given by%
\begin{equation}
U_{\mathfrak{M}}\left( s,u\right) =-sgn\left( u\kappa \right) B\left(
s\right) .  \tag{3.4}
\end{equation}%
Next the curvatures (mean and Gaussian) and the matrix of shape operator of
this surface are respectively as follows%
\begin{equation}
H_{\mathfrak{M}}=\frac{-sgn\left( u\kappa \right) \tau }{2u\kappa }\text{, \
\ and \ \ }K_{\mathfrak{M}}=0  \tag{3.5}
\end{equation}%
and%
\begin{equation}
S_{\mathfrak{M}}=\left( \frac{-sgn\left( u\kappa \right) \tau }{2u\kappa }%
\right) 
\begin{bmatrix}
1 & 1 \\ 
0 & 0%
\end{bmatrix}%
.  \tag{3.6}
\end{equation}%
Finally, the Christoffel symbols of the surface $\mathfrak{M}\left(
s,u\right) $ are given by%
\begin{eqnarray}
\left( \Gamma _{11}^{1}\right) _{\mathfrak{M}} &=&\dfrac{u\kappa _{s}+\kappa 
}{u\kappa },  \notag \\
\left( \Gamma _{11}^{2}\right) _{\mathfrak{M}} &=&\dfrac{-\kappa \left(
1+\left( u\kappa \right) ^{2}\right) -u\kappa _{s}}{u\kappa },  \TCItag{3.7}
\\
\left( \Gamma _{12}^{1}\right) _{\mathfrak{M}} &=&-\left( \Gamma
_{12}^{2}\right) _{\mathfrak{M}}=\dfrac{1}{u},  \notag \\
\left( \Gamma _{22}^{1}\right) _{\mathfrak{M}} &=&\left( \Gamma
_{22}^{2}\right) _{\mathfrak{M}}=0.  \notag
\end{eqnarray}

\bigskip

\section{The inverse surface of the tangent developable}

\bigskip

We show that $\mathfrak{N}$ is the inverse surface of the tangent
developable surface $\mathfrak{M}$ with respect to the inversion $\Phi .$
Thus the inverse surface $\mathfrak{N}$ has following parametrization 
\begin{equation}
\mathfrak{N}=c+\frac{r^{2}}{\left\Vert \mathfrak{M}-c\right\Vert ^{2}}\left( 
\mathfrak{M}-c\right) .  \tag{4.1}
\end{equation}

Hence, if we take into account the equalities $\left( 2.3\right) $ and $%
\left( 2.4\right) ,$ then the coefficients of the first and second
fundamental forms of the inverse surface $\mathfrak{N}$ by the help of these
of the surface $\mathfrak{M}$ are given by 
\begin{eqnarray}
E_{\mathfrak{N}} &=&\lambda ^{2}\left( 1+\left( u\kappa \right) ^{2}\right) ,%
\text{ }F_{\mathfrak{N}}=G_{\mathfrak{N}}=\lambda ^{2},  \TCItag{4.2} \\
l_{\mathfrak{N}} &=&sgn\left( u\kappa \right) \lambda u\kappa \tau -2\delta
\left( 1+\left( u\kappa \right) ^{2}\right) ,\text{ }m_{\mathfrak{N}}=n_{%
\mathfrak{N}}=-2\delta ,  \TCItag{4.3}
\end{eqnarray}%
where $E_{\mathfrak{N}},$ $F_{\mathfrak{N}},$ $G_{\mathfrak{N}}$ and $l_{%
\mathfrak{N}},$ $m_{\mathfrak{N}},$ $n_{\mathfrak{N}}$ are the coefficients
of the first and second fundamental forms of the inverse surface $\mathfrak{N%
},$ respectively.

Morever,the Gauss and mean curvatures of the inverse surface $\mathfrak{N}$
by the help of these of the surface $\mathfrak{M}$ are respectively, using
by $\left( 2.5\right) $ and $\left( 2.6\right) ,$%
\begin{eqnarray}
K_{\mathfrak{N}} &=&\frac{4}{r^{2}}\eta \left( -sgn\left( u\kappa \right) 
\frac{\tau }{2\lambda u\kappa }+\frac{\eta }{r^{2}}\right) ,  \TCItag{4.4} \\
H_{\mathfrak{N}} &=&sgn\left( u\kappa \right) \frac{\tau }{2\lambda u\kappa }%
-\frac{2\eta }{r^{2}},  \TCItag{4.5}
\end{eqnarray}%
where $K_{\mathfrak{N}}$ and $H_{\mathfrak{N}}$ are the Gauss and the mean
curvatures of the inverse surface $\mathfrak{N},$ respectively.\newline

\bigskip

\textbf{Theorem 4.1. }Let $\mathfrak{N}$ be the inverse surface of the
tangent developable surface $\mathfrak{M}$ with respect to the inversion $%
\Phi .$ Denote by $S_{\mathfrak{N}}$ the matrix of the shape operator of the
inverse surface $\mathfrak{N}$, then $S_{\mathfrak{N}}$ is given by the help
of that of $\mathfrak{M}$ as follows

\begin{equation}
S_{\mathfrak{N}}=%
\begin{bmatrix}
sgn\left( u\kappa \right) \frac{\tau }{\lambda u\kappa }-\frac{2\eta }{r^{2}}
& sgn\left( v\kappa \right) \frac{\tau }{\lambda v\kappa } \\ 
0 & -\frac{2\eta }{r^{2}}%
\end{bmatrix}%
.  \tag{4.6}
\end{equation}

\bigskip

\textbf{Proof. }Let $S_{\mathfrak{M}}$ be the matrix of the shape operator
of surface $\mathfrak{M.}$ By using the equalities $\left( 2.3\right) $ and $%
\left( 2.4\right) ,$ we can write%
\begin{equation}
S_{\mathfrak{N}}\circ \Phi _{\ast }=-\lambda ^{-1}S_{\mathfrak{M}}-\frac{2}{%
r^{2}}\eta I_{2},  \tag{4.7}
\end{equation}%
\newline
where $I_{2}$ is identity, $\lambda =\tfrac{r^{2}}{\left\Vert \mathfrak{M}%
-c\right\Vert ^{2}}$ and $\eta =\left\langle U_{\mathfrak{M}},\left( 
\mathfrak{M}-c\right) \right\rangle .$ Hence from $\left( 3.6\right) $ and $%
\left( 4.7\right) ,$ we obtain that the equality $\left( 4.6\right) $ is
satisfied.\newline

\bigskip

\textbf{Theorem 4.2. }Let $\left( \Gamma _{jk}^{i}\right) _{\mathfrak{N}}$
be the Christoffel symbols of the inverse surface $\mathfrak{N}.$ The
Christoffel symbols of the inverse surface $\mathfrak{N}$ by the help of
these of the surface $\mathfrak{M}$ are given by

\begin{eqnarray*}
\left( \Gamma _{11}^{1}\right) _{\mathfrak{N}} &=&\dfrac{u\kappa _{s}+\kappa 
}{u\kappa }+\frac{\left( \left( u\kappa \right) ^{2}-1\right) \frac{\partial
\lambda ^{2}}{\partial s}+\left( \left( u\kappa \right) ^{2}+1\right) \frac{%
\partial \lambda ^{2}}{\partial u}}{2\lambda ^{2}\left( u\kappa \right) ^{2}}%
\text{,} \\
\left( \Gamma _{11}^{2}\right) _{\mathfrak{N}} &=&\dfrac{-\kappa \left(
1+\left( u\kappa \right) ^{2}\right) -u\kappa _{s}}{u\kappa }+\frac{\left(
\left( u\kappa \right) ^{2}+1\right) \frac{\partial \lambda ^{2}}{\partial s}%
+\left( \left( u\kappa \right) ^{2}+1\right) ^{2}\frac{\partial \lambda ^{2}%
}{\partial u}}{2\lambda ^{2}\left( u\kappa \right) ^{2}}, \\
\text{\ }\left( \Gamma _{12}^{1}\right) _{\mathfrak{N}} &=&\dfrac{1}{u}+%
\frac{\left( \left( u\kappa \right) ^{2}+1\right) \frac{\partial \lambda ^{2}%
}{\partial u}-\frac{\partial \lambda ^{2}}{\partial s}}{2\lambda ^{2}\left(
u\kappa \right) ^{2}}, \\
\text{\ }\left( \Gamma _{12}^{2}\right) _{\mathfrak{N}} &=&-\dfrac{1}{u}+%
\frac{\left( \left( u\kappa \right) ^{2}+1\right) \left( \frac{\partial
\lambda ^{2}}{\partial s}-\frac{\partial \lambda ^{2}}{\partial u}\right) }{%
2\lambda ^{2}\left( u\kappa \right) ^{2}}, \\
\text{ \ }\left( \Gamma _{22}^{1}\right) _{\mathfrak{N}} &=&\frac{\left( 
\frac{\partial \lambda ^{2}}{\partial u}-\frac{\partial \lambda ^{2}}{%
\partial s}\right) }{2\lambda ^{2}\left( u\kappa \right) ^{2}}, \\
\left( \Gamma _{22}^{2}\right) _{\mathfrak{N}} &=&\frac{\left( \left(
u\kappa \right) ^{2}-1\right) \frac{\partial \lambda ^{2}}{\partial u}+\frac{%
\partial \lambda ^{2}}{\partial s}}{2\lambda ^{2}\left( u\kappa \right) ^{2}}%
.
\end{eqnarray*}

\bigskip

\textbf{Proof. }Considering the equality $\left( 2.3\right) ,$ for $i,j,k=1,$
we can write%
\begin{equation*}
\left( \Gamma _{11}^{1}\right) _{\mathfrak{N}}=\left( \Gamma
_{11}^{1}\right) _{\mathfrak{M}}+\dfrac{\left[ E_{\mathfrak{M}}G_{\mathfrak{M%
}}-2F_{\mathfrak{M}}^{2}\right] \dfrac{\partial \lambda ^{2}}{\partial s}+F_{%
\mathfrak{M}}E_{\mathfrak{M}}\dfrac{\partial \lambda ^{2}}{\partial u}}{%
2\lambda ^{2}\left[ E_{\mathfrak{M}}G_{\mathfrak{M}}-F_{\mathfrak{M}}^{2}%
\right] },
\end{equation*}%
where $\lambda =\tfrac{r^{2}}{\left\Vert \mathfrak{M}-c\right\Vert ^{2}}$
and $\left( \Gamma _{11}^{1}\right) _{\mathfrak{M}}$ is the Christoffel
symbol of the tangent developable surface. Thus, from the equalities $\left(
3.2\right) $ and $\left( 3.7\right) ,$ we obtain%
\begin{equation*}
\left( \Gamma _{11}^{1}\right) _{\mathfrak{N}}=\dfrac{u\kappa _{s}+\kappa }{%
u\kappa }+\frac{\left( \left( u\kappa \right) ^{2}-1\right) \frac{\partial
\lambda ^{2}}{\partial s}+\left( \left( u\kappa \right) ^{2}+1\right) \frac{%
\partial \lambda ^{2}}{\partial u}}{2\lambda ^{2}\left( u\kappa \right) ^{2}}
\end{equation*}

Others are found in similar way.\newline

\bigskip

\textbf{Theorem 4.3 }Let $\mathfrak{N}$ be the inverse surface of the
tangent developable surface $\mathfrak{M}$ with respect to the inversion $%
\Phi .$ Then the inverse surface $\mathfrak{N}$ is a flat surface if and
only if either the normal lines to the surface $\mathfrak{M}$ or the tangent
planes of the surface $\mathfrak{M}$ pass through the center of inversion. 
\newline

\bigskip

\textbf{Proof. }Let us assume that the inverse surface $\mathfrak{N}$ is
flat, then from $\left( 4.4\right) ,$ we can write%
\begin{equation}
\frac{4}{r^{2}}\eta \left( -sgn\left( u\kappa \right) \frac{\tau }{2\lambda
u\kappa }+\frac{\eta }{r^{2}}\right) =0,  \tag{4.9}
\end{equation}%
where either 
\begin{equation}
\eta =\left\langle U_{\mathfrak{M}},\left( \mathfrak{M}-c\right)
\right\rangle =0,  \tag{4.10}
\end{equation}%
or%
\begin{equation}
sgn\left( u\kappa \right) \frac{\tau }{2\lambda u\kappa }=\frac{\eta }{r^{2}}%
.  \tag{4.11}
\end{equation}%
If the equality $\left( 4.10\right) $ is satisfied, then the tangent planes
of the surface $\mathfrak{M}$ pass through the center of inversion. If the
equality $\left( 4.11\right) $ holds, then it follows 
\begin{equation*}
U_{\mathfrak{M}}=sgn\left( u\kappa \right) \frac{\tau }{2u\kappa }\left( 
\mathfrak{M}-c\right) .
\end{equation*}%
Namely, the normal lines to the surface $\mathfrak{M}$ pass through the
center of inversion.

The proof of sufficient condition is obvious.\newline

\bigskip

\textbf{Theorem 4.4 }Let $\mathfrak{N}$ be the inverse surface of the
tangent developable surface $\mathfrak{M}$ with respect to $S_{c}\left(
r\right) .$ The inverse surface $\mathfrak{N}$ is minimal if and only if the
normal lines to the surface $\mathfrak{M}$ pass through the center of
inversion \newline

\bigskip

\textbf{Proof. }The proof is same with that of Theorem 4.1.

\bigskip

\textbf{Applications. }

\bigskip

\begin{gather*}
\FRAME{itbpF}{1.8663in}{1.7858in}{0in}{}{}{Figure}{\special{language
"Scientific Word";type "GRAPHIC";maintain-aspect-ratio TRUE;display
"USEDEF";valid_file "T";width 1.8663in;height 1.7858in;depth
0in;original-width 4.0352in;original-height 3.8588in;cropleft "0";croptop
"1";cropright "1";cropbottom "0";tempfilename
'LZUH4000.wmf';tempfile-properties "XPR";}} \\
\text{Fig }1.\text{ }
\end{gather*}%
\begin{gather*}
\FRAME{itbpF}{1.8308in}{1.8741in}{0in}{}{}{Figure}{\special{language
"Scientific Word";type "GRAPHIC";maintain-aspect-ratio TRUE;display
"USEDEF";valid_file "T";width 1.8308in;height 1.8741in;depth
0in;original-width 3.7299in;original-height 3.8181in;cropleft "0";croptop
"1";cropright "1";cropbottom "0";tempfilename
'LZUH5101.wmf';tempfile-properties "XPR";}} \\
\text{Fig }2.\text{ }
\end{gather*}

\bigskip

\textbf{Fig 1: }The helicoid given by $\left( u\cos v,u\sin v,2v\right) .$

\bigskip

\textbf{Fig 2: }The inverse surface of the helicoid with respect to unit
sphere given by 
\begin{equation*}
\left( \frac{u}{u^{2}+4v^{2}}\cos v,\frac{u}{u^{2}+4v^{2}}\sin v,\frac{2v}{%
u^{2}+4v^{2}}\right) .
\end{equation*}

\end{document}